\documentclass[a4paper,11pt]{article}

\usepackage[utf8]{inputenc}
\usepackage[left=2cm,right=2cm]{geometry}
\usepackage{mathtools}
\usepackage{amssymb}
\usepackage{amsthm}
\usepackage{hyperref}
\usepackage{doi}
\usepackage[numbers,sort&compress]{natbib}
\usepackage{color}

\hypersetup{hidelinks}

\newtheorem{theorem}{Theorem}

\title{\bf Uniqueness in an incompressible fluid with objective heat conduction model}
\author{Kinza Arshad, Martina Nunziata, Vincenzo Tibullo\footnote{%
    Corresponding author: Vincenzo Tibullo, \texttt{vtibullo@unisa.it}}\\
\small{Dipartimento di Matematica, Università di Salerno, Italy}}
\date{}

\begin{document}

\maketitle

\begin{abstract}
In the context of heat conduction, the equations of Maxwell-Cattaneo is a successful alternative to the Fourier's law. In 2009, 
%\citet{Christov2009} 
Christov
proposed a frame-independent (or objective) generalization of Maxwell-Cattaneo equation. For a more general class of objective heat conduction Maxwell-Cattaneo-Christov equations proposed by 
%\citet{Angeles2023}, 
Angeles,
a uniqueness result is proved, applied to an incompressible fluid.
\end{abstract}

\section{Introduction}
Over the past few decades, researchers have found great interest in solving the heat conduction problem. In 1822, \citet{Fourier1888} introduced the law of heat conduction, which states that ``heat transfer in a medium occurs at an inertial rate''. This law along with the energy balance equation leads to a parabolic system of equations for heat conduction. However, a key issue is that no real material or body perfectly conforms to Fourier's law of heat conduction. Many researches have been conducted to verify that Fourier's law is not applicable in certain circumstances (see e.g \citep{Yang2010}). Several alternative models have been developed to correct its unrealistic features, among which the Cattaneo-Maxwell heat transfer law \cite{Joseph1989} is one of the most recognized 
\[
\textcolor{black}{ \tau \dot{q}_i + q_i = - \kappa T_{,i}, }
\]
where \textcolor{black}{\(q_i\) are the components of} the heat flux, \textcolor{black}{ \(\dot{q}_i=\partial q_i/\partial t\) }denotes the partial time derivative of heat flux,  \(\kappa \) is the thermal conductivity and \(T\) stands for temperature,  \(\tau >0 \) is a constant and serves as an intrinsic relaxation time, representing the time delay necessary for heat conduction to occur within a volume element after a temperature gradient has been established. Consequently, this new term introduces a form of \textquotedblleft thermal inertia\textquotedblright.  In this modified version of Fourier's law, heat transfer in the medium is not instantaneous; rather, it occurs through the propagation of thermal waves at a finite speed, a phenomenon known as second sound \cite{Jordan2014, Straughan2011}.

According to \citet{Christov2005}, the above equation  contradicts a fundamental principle of classical mechanics, leading to paradoxical interpretations of the propagation of thermal waves. This inconsistency arises from the notion that thermal inertia should be an inherent property of the material point. Consequently, the rate of change of heat flux over time must consider both the change occurring at the geometric point (the partial time derivative) and the transport of material quantities when the medium is in motion. To resolve this problem, they proposed replacing the partial time derivative with a material derivative. This perspective enabled \citet{Christov2009} to formulate a frame-independent version of the Cattaneo-Maxwell law, which substitutes the partial time derivative of heat flux with a Lie-Oldroyd upper convected material (or objective) derivative \cite{Oldroyd1950}. The resulting formulation is expressed as follows
\[
\textcolor{black}{ \tau \left[ \dot{q}_i + v_j q_{i,j} - v_{i,j} q_j + v_{j,j} q_i \right] + q_i = -\kappa T_{,i}. }
\]
This formulation satisfies the principle of objectivity in continuum mechanics. 
However, it has been recently questioned that this system is not hyperbolic \cite{Angeles2022}. 
Hyperbolicity in quasilinear systems refers to the possibility of identifying multiple linearly independent waves traveling in any given spatial direction, and for the Cattaneo-Christov system there are particular directions where these distinct waves do not exist. 
Notably, since the characteristic speeds of this system are real, it is categorized as weakly hyperbolic.  
While the Cauchy problem \cite{Colombini2018} for weakly hyperbolic systems of equations can be well-posed in Gevrey spaces, this condition alone is not sufficient for overall well-posedness. 
\textcolor{black}{ For the well-posedness of weakly hyperbolic Cauchy problems, see also \citet{Colombini2018} and the references therein.}
Additionally, it is well established that hyperbolicity in first-order quasilinear systems is a necessary requirement for the existence of energy estimates \cite{Dafermos2010, Serre1999}. 
\textcolor{black}{Although systems that are only weakly hyperbolic can have smooth solutions (see Theorem 5, p. 430 of \citet{Evans2022}), they cannot accept \(L^2\) energy estimates (see Theorem 2.2.3 in \citet{Kreiss2004}). 
In fact, the \(L^2\) well-posedness requires the hyperbolicity property (see \citet{Kajitani1973}).}
As a result of these considerations \citet{Angeles2023} proposed the evolution equation 
\begin{align}
    \label{eq:1}
    \tau \left(\dot{q}_i + v_j q_{i,j}-\frac{1}{2}(v_{i,j}-v_{j,i})q_j+\frac{\lambda}{2}(v_{i,j}+v_{j,i})q_j+ \mu v_{j,j}q_i\right) + q_i= -\kappa T_{,i}.
\end{align}
\citet{Morro2018} reveals that for any pair of invariant scalars \((\lambda, \mu) \) the bracketed term in \eqref{eq:1} constitutes an objective derivative. Furthermore, when \((\lambda, \mu)\) are considered as real-valued functions of \((\rho, T)\), the quasilinear coupling equation model \eqref{eq:1} with conservation of mass, the balance of linear momentum and the balance of total energy (for a compatible inviscid fluid) % rephrase this again .
is a hyperbolic system of equations if and only if \((\lambda(\rho, T), \mu(\rho, T)) = (1, -1)\) for every \(\rho > 0\) and \(T > 0\) \textcolor{black}{ \cite{Angeles2023}. Moreover, in \citet{Angeles2024} it is shown that the \((1, -1)\) quasilinear model is the only one with a symbolic symmetrizer and thus, only this system has the possibility of have a unique well-posed solution in \(L^2\) and Sobolev spaces. }

Our aim in this article is to carry forward the result attained by \citet{Tibullo2011} in the framework of objective heat conduction models and investigate the corresponding uniqueness result. Meanwhile, we want to scrutinize this model when subjected to incompressibility condition. This can be achieved by proving the uniqueness of solution for the given model.

\section{Basic equations}
We start from a model of a fluid in a bounded region \(\Omega\) 
that satisfy the following equations (see \cite{Tibullo2011})
\begin{enumerate}
    \item the balance of linear momentum
        \begin{equation}
        \label{eq:linear momentum}
         v_{i,t}+v_jv_{i,j}=-\frac{1}{\rho}p_{,i}+\nu v_{i,jj}+\alpha gk_iT,
        \end{equation}
        where \(v_i\) is the velocity field, \(\rho\) is the mass density, \(p\) is the pressure, \(\alpha\) is the thermal expansion coefficient, \(g\) is acceleration due to gravity, \(\nu\) is the kinematic viscosity, \(k_i\) is the unit vector of the vertical \(z\)-axis and \(T\) is the absolute temperature;
    \item the balance of energy
        \begin{equation}
         \label{eq:energy}
        \rho c_P(T_{,t}+v_{i}T_{,i})=-q_{i,i},
        \end{equation}
        where \(c_P\) is the specific heat at constant pressure and \(q_i\) is the heat flux;
    \item the fluid is supposed to be incompressible, so that
        \begin{equation}
        \label{eq:incompressibility}
         v_{i,i}=0.
        \end{equation}
        Moreover, we suppose that the density \(\rho\) is uniform in the reference configuration, so it remains uniform because of incompressibility.
        \end{enumerate}
For the heat conduction model, we chose that one proposed by \citet{Angeles2023}
\begin{equation*}
    \label{eq:MCC objective heat conduction-1}
    \tau\left(q_{i,t}+v_jq_{i,j}-\frac{1}{2}(v_{i,j}-v_{j,i})q_j+\frac{\lambda}{2}(v_{i,j}+v_{j,i})q_j+\mu v_{j,j}q_i\right)+q_i=-\kappa T_{,i},
\end{equation*}
where \(\tau\) is the thermal relaxation time, \(\kappa\) is the thermal conductivity.% and \((\lambda,\mu)\) is a pair of parameters and it is proved that this heat conduction model is objective for every choice of the pair. In the article \cite{Angeles2023} is also written that only when \((\lambda,\mu)=(1,-1)\) yields a hyperbolic system of quasilinear equations.

When the fluid is incompressible, the parameter \(\mu\) disappears and the objective heat conduction model is determine by the values of \(\lambda\)
\begin{equation}
    \label{eq:MCC objective heat conduction}
    \tau\left(q_{i,t}+v_jq_{i,j}-\frac{1}{2}(v_{i,j}-v_{j,i})q_j+\frac{\lambda}{2}(v_{i,j}+v_{j,i})q_j\right)+q_i=-\kappa T_{,i}.
\end{equation}

Our aim is to extend the results obtained in \cite{Tibullo2011} in the context of objective heat conduction and to analyze the corresponding uniqueness result. \textcolor{black}{ To this aim, we had to overcome some technical difficulties due to the more complicated form of the equations. We remark that we maintain an arbitrary value of \(\lambda\) to show that the results obtained are independent of its value.}

\section{Uniqueness result}

We can prove this uniqueness result associated to the model described above.
\textcolor{black}{By a regular solution we mean a set of functions \((v_i, T, q_i, p)\) such that \(v_i\in C^{2,1}(\Omega\times I)\), \(T, q_i \in C^{1,1}(\Omega\times I)\) and \(p\in C^{1,0}(\Omega\times I)\), where \(I\) is the time interval in which we seek for a solution.}
\begin{theorem}
    For an incompressible fluid and in the context of the objective heat conduction model represented by eqs.\eqref{eq:linear momentum}-\eqref{eq:MCC objective heat conduction}, there exists at most one regular solution \((v_i, T, q_i, p)\) for the initial-boundary value problem with assigned and arbitrary boundary and initial conditions in a bounded domain \(\Omega\).
\end{theorem}
\begin{proof} 
Let us consider two solutions of the problem
\[(v_i^{(1)}, T^{(1)}, q_i^{(1)}, p^{(1)}),\quad (v_i^{(2)}, T^{(2)}, q_i^{(2)}, p^{(2)})\]
that satisfy the same boundary and initial conditions. We define
\begin{equation}
    \label{eq:differences of solutions}
    \begin{aligned}
        u_i&=v_i^{(1)}-v_i^{(2)},\\
        \theta&=T^{(1)}-T^{(2)},\\
        Q_i&=q_i^{(1)}- q_i^{(2)},\\
        \pi&=p^{(1)}-p^{(2)},
    \end{aligned}
\end{equation}
so we have
\begin{equation}
    \label{eq: BC and IC}
    u_i=0,\quad \theta=0,\quad Q_i=0,\quad \pi=0,
\end{equation}
    on \(\partial\Omega\times[0,+\infty)\) and on \(\Omega\times\{0\}\). It is easy to see that \((u_i, \theta, Q_i, \pi)\) satisfies the following system of equations
    \begin{equation}
        \label{eq:new-system}
        \begin{aligned}
            & u_{i,t}+u_jv_{i,j}^{(1)}+v_j^{(2)}u_{i,j}=-\frac{1}{\rho}\pi_{,i}+\nu u_{i,jj}+\alpha gk_i\theta,\\
            & \rho c_P(\theta_{,t}+u_iT_{,i}^{(1)}+v_{i}^{(2)}\theta_{,i})=-Q_{i,i},\\
            & \tau\bigg[Q_{i,t}+u_jq_{i,j}^{(1)}+v_j^{(2)}Q_{i,j}\\
            & \qquad +\bigg(-\frac{1}{2}(u_{i,j}-u_{j,i})+\frac{\lambda}{2}(u_{i,j}+u_{j,i})\bigg)q_j^{(1)}\\
            & \qquad +\bigg(-\frac{1}{2}(v_{i,j}^{(2)}-v_{j,i}^{(2)})+\frac{\lambda}{2}(v_{i,j}^{(2)}+u_{j,i}^{(2)})\bigg)Q_j\bigg]+Q_i=-\kappa\theta_{,i}.
        \end{aligned}
    \end{equation}
Let \(\theta_0\) be the temperature in the reference configuration and multiply the equations of system \eqref{eq:new-system} by \(\rho u_i, \theta/\theta_0\) and \(Q_i/\kappa\theta_0\), respectively. Next we sum over the repeated index and integrate over \(\Omega\). Thanks to the incompresibility and the boundary conditions \eqref{eq: BC and IC} we have
\[
\begin{aligned}
    &\frac{1}{2}\rho\frac{d}{dt}\|\mathbf{u}\|^2+\rho\int_{\Omega}u_iu_jv_{i,j}^{(1)}\,dx+\rho\int_{\Omega}u_iv_j^{(2)}u_{i,j}\,dx=-\rho\nu\|\nabla\mathbf{u}\|^2+\rho \alpha g\int_{\Omega}u_ik_i\theta\,dx,\\
    &\frac{1}{2}\frac{\rho c_P}{\theta_0}\frac{d}{dt}\|\theta\|^2+\frac{\rho c_P}{\theta_0}\int_{\Omega} \theta u_i T_{,i}^{(1)}\,dx+\frac{\rho c_P}{\theta_0}\int_{\Omega} \theta v_i^{(2)} \theta_{,i}\,dx=-\frac{1}{\theta_0}\int_{\Omega}\theta Q_{i,i}\,dx,\\
    &\frac{\tau}{\kappa\theta_0}\frac{1}{2}\frac{d}{dt}\|\mathbf{Q}\|^2+\frac{\tau}{\kappa\theta_0}\int_{\Omega}Q_iu_jq_{i,j}^{(1)}\,dx+\frac{\tau}{\kappa\theta_0}\int_{\Omega}Q_iv_j^{(2)}Q_{i,j}\,dx\\
    &\qquad+\frac{\tau}{\kappa\theta_0}\int_{\Omega}Q_i\bigg(-\frac{1}{2}(v_{i,j}^{(2)}-v_{j,i}^{(2)})+\frac{\lambda}{2}(v_{i,j}^{(2)}+v_{j,i}^{(2)})\bigg)Q_j\,dx\\
    &\qquad+\frac{\tau}{\kappa\theta_0}\int_{\Omega}Q_i\bigg(-\frac{1}{2}(u_{i,j}-u_{j,i})+\frac{\lambda}{2}(u_{i,j}+u_{j,i})\bigg)q_j^{(1)}\,dx+\frac{1}{\kappa\theta_0}\|\mathbf{Q}\|^2=-\frac{1}{\theta_0}\int_{\Omega}Q_i\theta_{,i}\,dx,
\end{aligned}
\]
where the norms of \(\mathbf{u}, \nabla\mathbf{u}, \theta\) and \(\mathbf{Q}\) are defined in the \(L^2\) spaces as follows
\[\begin{aligned}
    &\|\mathbf{u}\|^2=\int_{\Omega} u_iu_i\,dx,\qquad \|\nabla\mathbf{u}\|^2=\int_{\Omega} u_{i,j}u_{i,j}\,dx,\\
    &\|\theta\|^2=\int_{\Omega} \theta^2\,dx,\quad\qquad \|\mathbf{Q}\|^2=\int_{\Omega} Q_iQ_i\,dx.
\end{aligned}\]
Adding together the above equations, we have
\begin{equation}
    \label{eq:together norm-1}
    \begin{aligned}
    \frac{1}{2}\frac{d}{dt}&\left(\rho\|\mathbf{u}\|^2+\frac{\rho c_P}{\theta_0}\|\theta\|^2+\frac{\tau}{\kappa\theta_0}\|\mathbf{Q}\|^2\right)=\\
    &-\rho\int_{\Omega}u_iu_jv_{i,j}^{(1)}\,dx-\rho\nu\|\nabla\mathbf{u}\|^2+\rho \alpha g\int_{\Omega}u_ik_i\theta\,dx\\
    &-\frac{\rho c_P}{\theta_0}\int_{\Omega} \theta u_i T_{,i}^{(1)}\,dx-\frac{\tau}{\kappa\theta_0}\int_{\Omega}Q_iu_jq_{i,j}^{(1)}\,dx-\frac{1}{\kappa\theta_0}\|\mathbf{Q}\|^2\\
    &-\frac{\tau}{\kappa\theta_0}\int_{\Omega}Q_i\bigg(-\frac{1}{2}(v_{i,j}^{(2)}-v_{j,i}^{(2)})+\frac{\lambda}{2}(v_{i,j}^{(2)}+v_{j,i}^{(2)})\bigg)Q_j\,dx\\
    &-\frac{\tau}{\kappa\theta_0}\int_{\Omega}Q_i\bigg(-\frac{1}{2}(u_{i,j}-u_{j,i})+\frac{\lambda}{2}(u_{i,j}+u_{j,i})\bigg)q_j^{(1)}\,dx,
    \end{aligned}
\end{equation}
where we use, on some therms, the divergence theorem, the boundary conditions \eqref{eq: BC and IC} and the incompressibility. 

We can define \(V_1, V_2, \Theta, Q, Q'\) the finite positive maxima of the functions 
\[|v_{i,j}^{(1)}|,\quad |v_{i,j}^{(2)}|,\quad |\theta_{,i}^{(1)}|,\quad |q_{j}^{(1)}|,\quad |q_{i,j}^{(1)}|,\qquad \forall\mathbf{x}\in\Omega,\,\forall t\in[0,t_{\max}],\]
since they are continuous functions on a compact set\textcolor{black}{; here \(t_{\max}\) is an arbitrary time for which the solution exists}. To perform some estimates on the terms of eq. \eqref{eq:together norm-1}, we can apply the arithmetic-geometric mean inequality, i.e.
\[a_ib_i\leq\frac{1}{2}\left(\frac{1}{\epsilon}a_ia_i+\epsilon b_ib_i\right),\quad\forall\epsilon>0,\]
where \(a_i\) and \(b_i\) are arbitrary vectors. We can show that
\begin{equation*}
    \begin{aligned}
        &-\rho\int_{\Omega}u_iu_jv_{i,j}^{(1)}\,dx\leq3\rho V_1\|\mathbf{u}\|^2\,,\\
        &+\rho \alpha g\int_{\Omega}u_ik_i\theta\,dx\leq\frac{1}{2}\rho\alpha g\left(\textcolor{black}{\|\mathbf{u}\|^2+\|\theta\|^2}\right)\,,\\
        &-\frac{\rho c_P}{\theta_0}\int_{\Omega} \theta u_i T_{,i}^{(1)}\,dx\leq\frac{\rho c_P\Theta}{2\theta_0}\left(\textcolor{black}{ 3\|\theta\|^2 +\|\mathbf{u}\|^2 }\right)\,,\\
        &-\frac{\tau}{\kappa\theta_0}\int_{\Omega}Q_iu_jq_{i,j}^{(1)}\,dx\leq\frac{3\tau Q'}{2\kappa\theta_0}\left(\textcolor{black}{ \|\mathbf{Q}\|^2+\|\mathbf{u}\|^2 }\right)\,,\\
        &-\frac{\tau}{\kappa\theta_0}\int_{\Omega}Q_i\bigg(-\frac{1}{2}(v_{i,j}^{(2)}-v_{j,i}^{(2)})+\frac{\lambda}{2}(v_{i,j}^{(2)}+v_{j,i}^{(2)})\bigg)Q_j\,dx\leq\frac{3\tau V_2\lambda}{\kappa\theta_0}\|\mathbf{Q}\|^2\,,\\
        &-\frac{\tau}{\kappa\theta_0}\int_{\Omega}Q_i\bigg(-\frac{1}{2}(u_{i,j}-u_{j,i})+\frac{\lambda}{2}(u_{i,j}+u_{j,i})\bigg)q_j^{(1)}\,dx\leq\\
        &\qquad+\frac{\tau Q}{\kappa\theta_0}\frac{\lambda-1}{4}\left(\frac{3}{\textcolor{black}{\epsilon_1}}\|\mathbf{Q}\|^2+\textcolor{black}{\epsilon_1}\|\nabla\mathbf{u}\|^2\right)+\frac{\tau Q}{\kappa\theta_0}\frac{\lambda+1}{4}\left(\frac{3}{\textcolor{black}{\epsilon_2}}\|\mathbf{Q}\|^2+\textcolor{black}{\epsilon_2}\|\nabla\mathbf{u}\|^2\right)\,.
    \end{aligned}
\end{equation*}
Using these estimates, \textcolor{black}{ one takes the values \(\epsilon_1\) and \(\epsilon_2\) so that the right hand side terms in the estimate can be absorbed in the left hand side of \eqref{eq:together norm-1} with \(\rho\nu\|\nabla\mathbf{u}\|^2\), and }the estimate will be
\begin{equation}
    \label{eq:estimate ABC-1}
    \frac{1}{2}\frac{d}{dt}\left(\rho\|\mathbf{u}\|^2+\frac{\rho c_P}{\theta_0}\|\theta\|^2+\frac{\tau}{\kappa\theta_0}\|\mathbf{Q}\|^2\right)\leq A\|\mathbf{u}\|^2+B\|\theta\|^2+C\|\mathbf{Q}\|^2\,,
\end{equation}
where the exact positive values of \(A, B, C\) can be deduced from the above estimates. If we define 
\begin{equation}
    \label{eq:definition F-1}
    F(t)=\frac{1}{2}\left(\rho\|\mathbf{u}\|^2+\frac{\rho c_P}{\theta_0}\|\theta\|^2+\frac{\tau}{\kappa\theta_0}\|\mathbf{Q}\|^2\right)
\end{equation}
we can find a positive constant \(\tau_0\) such that the inequality \eqref{eq:estimate ABC-1} is equivalent to
\begin{equation*}
    \dot{F}(t)\leq\frac{F(t)}{\tau_0}.
\end{equation*}
Since \(F(t)\geq0\) for all \(t\geq0\) and thanks to Gronwall's inequality applied to \(F(t)\), we find
\[F(t)\leq F(0)e^{-t/\tau_0},\qquad\forall t\geq0,\]
where \(F(0)=0\), so \(F(t)\) vanishes identically. The same happens to \(u_i\), \(\theta\) and \(Q_i\) and for \(\pi\), since we note from \eqref{eq:new-system}\(_1\) that \(\pi_{,i}=0\) so \(\pi=\text{const}\) and thanks to the boundary conditions \eqref{eq: BC and IC} we have \(\pi=0\) and we find that the two solutions \((v_i^{(1)}, T^{(1)}, q_i^{(1)}, p^{(1)}), (v_i^{(2)}, T^{(2)}, q_i^{(2)}, p^{(2)})\) coincide everywhere in \(\Omega\).
\end{proof}

\section{Conclusions}
In this work we proved a uniqueness result for an incompressibile fluid in which we adopt an objective heat conduction model of Maxwell-Cattaneo-Christov type proposed by \citet{Angeles2023}. This result extend that one proved in \cite{Tibullo2011} for the choice of the parameter \(\lambda\) that defines a class of objective heat conduction models. For the future, our aim is to understand more deeply the class of objective derivatives and give them a meaning in relation to thermodynamic consistency.  

\section*{Acknowledgments}

This paper was partially supported by INdAM (GNFM), by the project  \textit{Non linear models for magma transport and volcanoes}, project code P20222B5P9, CUP:D53D23018760001, PRIN 2022 PNRR and by project \textit{Transport phonema in low dimensional structures: models, simulations and
theoretical aspects}, project code 2022TMW2PY, CUP:D53D23006000006, PRIN 2022.

\bibliography{references}
\bibliographystyle{unsrtnat}

\end{document}